\documentclass[a4paper,12pt]{article}
\usepackage[T1]{fontenc}
\usepackage[utf8]{inputenc}
\usepackage{lmodern}
\usepackage[margin=25mm]{geometry}
\usepackage{amsmath,amssymb,amsthm}
\usepackage{graphicx,booktabs,array,xcolor}
\usepackage{placeins,etoolbox}
\usepackage{xurl}
\usepackage[hidelinks]{hyperref}
\hypersetup{pdftitle={CA2 numbers and the Riemann hypothesis: a numerical investigation of Ramanujan's inequalities},pdfauthor={Dmitry R. Leshchiner, Oleg R. Musin, James Thompson}}

\newtheorem{theorem}{Theorem}
\theoremstyle{definition}
\newtheorem{question}{Question}
\newtheorem{conjecture}{Conjecture}
\newcommand{\KeepStatementTogether}[1]{%
  \BeforeBeginEnvironment{#1}{\par\begin{samepage}}%
  \AfterEndEnvironment{#1}{\end{samepage}}}
\forcsvlist{\KeepStatementTogether}{theorem,question,conjecture}
\numberwithin{equation}{section}
\newcommand{\CA}{\mathrm{CA}}
\newcommand{\RH}{\mathrm{RH}}
\newcommand{\eg}{e^{\gamma}}
\title{CA2 numbers and the Riemann hypothesis:\\a numerical investigation of\\Ramanujan's inequalities}
\author{Dmitry R. Leshchiner \and Oleg R. Musin \and James Thompson}
\date{}
\begin{document}
\maketitle

\begin{abstract}
Robin's criterion reduces the Riemann hypothesis to an inequality for colossally abundant (CA) numbers, and Wu's theorem reduces the test further to the subclass CA2. We study this reduction numerically and consider smaller subsets on which the same criterion remains sufficient. We also show that the logarithmic decay law suggested by Briggs's computations follows from Ramanujan's conditional bounds. To investigate the next term, we compare the normalized Robin difference with Ramanujan's oscillatory expression involving the nontrivial zeros of the zeta function. The resulting difference has a slowly varying component on both CA and CA2. Its scaled values continue to vary over the computed range, leading to questions about convergence, upper bounds, and higher-order terms. We describe the sampling used in the figures and estimate the error from truncating the sum over zeros.

\end{abstract}

\smallskip


\section{Introduction}
\label{sec:intro}
Let $\sigma(n)=\sum_{d\mid n}d$, and let $\gamma$ be the Euler--Mascheroni constant. Robin \cite{Robin} proved that the Riemann hypothesis (RH) is equivalent to
\begin{equation}\label{eq:robin}
 \frac{\sigma(n)}n<\eg\log\log n\qquad(n>5040).
\end{equation}
It is enough to test this inequality on the colossally abundant (CA) numbers, which maximize $\sigma(n)/n^{1+\epsilon}$ for some $\epsilon>0$. Wu \cite{Wu19} reduced the test to CA2, a subclass containing fewer than one in a thousand CA numbers in his reported computation. We examine the behavior of Robin's quotient on these classes and consider further reductions.

Briggs \cite{Briggs} verified Robin's inequality for $5040<n<10^{10^{10}}$ and investigated the distance from equality. Morrill and Platt \cite{MP18} subsequently extended the verified range to $10^{10^{13}}$. The logarithmic law suggested by Briggs's computations follows from Ramanujan's estimates under RH. These estimates occur in the continuation of Ramanujan's manuscript omitted from the 1915 publication; this part appeared in facsimile in 1988 and in an annotated journal edition in 1997 \cite{Ram97,NS}.

Ramanujan's work \cite{Ram15,Ram97} suggests the following normalization of $\sigma(n)/n$:
\begin{equation}\label{eq:T}
T(n)=\left(\eg\log\log n-\frac{\sigma(n)}n\right)\sqrt{\log n}.
\end{equation}
Under RH, the lower and upper limits of $T(n)$ along the CA numbers lie between two positive constants approximately equal to $1.3932$ and $1.5578$. Robin's criterion requires only $T(n)>0$ for every $n>5040$. It is also enough to prove this inequality for all sufficiently large CA numbers: if RH is false, Robin's oscillation theorem and the reduction to CA numbers give arbitrarily large CA counterexamples \cite{Robin,Musin}.

For a more precise comparison, we use the function $R(n)$ defined by a sum over zeta zeros in \eqref{eq:R}. Under RH, it lies between the same two constants. Our main numerical quantities are $D(n)=R(n)-T(n)$ and
\begin{equation}\label{eq:Q-intro}
 Q(n)=e^{-\gamma}\bigl(R(n)-T(n)\bigr)\log\log n.
\end{equation}
Ramanujan's estimate and its effective refinement by Nicolas \cite{Nic21} suggest this scale for the next correction. The figures show similar oscillations in $T$ and $R$, together with a more slowly varying difference. The observed increase in the scaled difference may be part of a longer oscillation; the available range does not determine a limiting coefficient.

A separate question is whether one can reduce CA2 to smaller sets on which Robin's criterion remains sufficient. Record values of $G(n)=\sigma(n)/(n\log\log n)$ give one such reduction. A different construction of nested subsets of CA numbers, using changes of coordinates and supporting lines, is developed in \cite{Musin26}. We discuss its relation to the CA2 problem in Section~\ref{sec:thinning}.

Section~\ref{sec:framework} gives the normalizations, the deduction of Briggs's logarithmic law, and Wu's reduction. Section~\ref{sec:experiments} presents the numerical observations and the questions they suggest. Section~\ref{sec:computation} describes generation, sampling, and truncation errors. Section~\ref{sec:asymptotic} discusses higher-order corrections and the sharpness of Ramanujan's constants.

\section*{Acknowledgments}
The authors thank Aiya Kuchukova for her contribution to the early stages of this project. ChatGPT (OpenAI) was used to assist with the editing of the manuscript.

\section{The mathematical framework}
\label{sec:framework}
\subsection{The Robin difference and the sum over zeta zeros}
All logarithms in the text are natural. Put
\[
 u=\log\log n,\qquad
 G(n)=\frac{\sigma(n)}{n\log\log n},\qquad n\ge3.
\]
Then $T(n)=(\eg-G(n))u e^{u/2}$. This identity explains why a plot of $G$ alone can conceal the variation we seek: differences of order one in $T$ correspond to differences of order $e^{-u/2}/u$ in $G$.

Under RH, write the nontrivial zeta zeros with positive imaginary part as $\rho_j=\tfrac12+it_j$, counted with multiplicity. The comparison function is
\begin{equation}\label{eq:R}
 R(n)=\eg\left(2(\sqrt2-1)
       -2\sum_{j\ge1}\frac{\cos(t_j u)}{\tfrac14+t_j^2}\right).
\end{equation}
The sum converges absolutely. With the common notation
\[
 S(x)=\sum_\rho\frac{x^{\rho-1}}{\rho(\rho-1)},
\]
formula~\eqref{eq:R} is $R(n)=\eg\bigl(2(\sqrt2-1)+\sqrt{\log n}\,S(\log n)\bigr)$, with conjugate zeros paired. To compare signs with Nicolas \cite{Nic21}, his notation is
\[
 S_1(x)=\sum_\rho\frac{x^{\rho-1}}{\rho(1-\rho)}=-S(x).
\]
Thus $R(n)=\eg\bigl(2(\sqrt2-1)-\sqrt{\log n}\,S_1(\log n)\bigr)$. Under RH, $\rho(1-\rho)=\tfrac14+t^2$; pairing conjugate zeros gives the minus sign in \eqref{eq:R}.

We use three different residual scales:
\begin{equation}\label{eq:residuals}
 D(n)=R(n)-T(n),\qquad P(n)=e^{-\gamma}D(n),\qquad Q(n)=uP(n).
\end{equation}
Figures labeled $R-T$ show $D$; a figure labeled $(R-T)\log\log n$ shows $\eg Q$. This constant factor matters when discussing a numerical bound for $Q$.

Set
\[
 \tau=2+\gamma-\log(4\pi)
      =2\sum_{j\ge1}\frac1{\tfrac14+t_j^2}
      =0.0461914179\ldots,
\]
where the zero-sum identity is written under RH. The resulting bounds for $R$ are
\begin{equation}\label{eq:constants}
 \begin{split}
 c_1&=\eg\bigl(2(\sqrt2-1)-\tau\bigr)=1.3932184417\ldots,\\
 c_2&=\eg\bigl(2(\sqrt2-1)+\tau\bigr)=1.5577589626\ldots,
 \end{split}
 \qquad c_1\le R(n)\le c_2.
\end{equation}
Ramanujan's conditional inequalities \cite{Ram97,Musin} give
\begin{equation}\label{eq:ramanujan}
 c_1\le R_L:=\liminf_{\substack{n\to\infty\\n\in\CA}}T(n)
 \le R_U:=\limsup_{\substack{n\to\infty\\n\in\CA}}T(n)\le c_2.
\end{equation}
These are asymptotic bounds; an individual computed value of $T$ below $c_1$ does not contradict them. The lower bound also holds for the lower limit over all integers, by the strong Ramanujan theorem \cite{Musin}. The upper bound is specific to suitable abundant sequences.

Ramanujan's estimate, made effective by Nicolas \cite{Nic21}, gives the one-sided asymptotic comparison
\begin{equation}\label{eq:one-sided}
 T(n)\ge R(n)-O\!\left(\frac1{\log\log n}\right)
 \qquad(n\to\infty,\ \RH).
\end{equation}
Thus $Q$ is bounded above for sufficiently large $n$ under RH. This one-sided estimate does not imply convergence or a lower bound for $Q$.

\subsection{Briggs's experiment and Ramanujan's inequalities}
\label{sec:briggs}
Put $\Delta(n)=\eg\log\log n-\sigma(n)/n = T(n)/\sqrt{\log n}$. The logarithmic law suggested by Briggs's computations \cite{Briggs} has a leading slope $-1/2$:
\begin{equation}\label{eq:briggs-law}
 \log\Delta(n)=-\frac12\log\log n+o(\log\log n)
 \qquad(n\to\infty,\ n\in\CA,\ \RH).
\end{equation}
In fact, Ramanujan's bounds give a stronger bounded-remainder statement
\begin{equation}\label{eq:briggs-strong}
 \log\Delta(n)=-\frac12\log\log n+O(1).
\end{equation}
Indeed, \eqref{eq:ramanujan} implies that $T(n)$ is eventually bounded above and bounded away from zero on the CA numbers. Since $\Delta(n)=T(n)/\sqrt{\log n}$, taking logarithms gives
\[
 \log\Delta(n)=-\frac12\log\log n+\log T(n),
\]
and $\log T(n)=O(1)$.

Thus Briggs's observed slope follows from Ramanujan's conditional bounds. The normalization $T$ removes this leading scale; comparison with $R$ then concerns the next term.

\subsection{CA numbers and Wu's reduction}
An integer $N$ is called \emph{colossally abundant} if, for some $\epsilon>0$,
\[
 \frac{\sigma(k)}{k^{1+\epsilon}}
 \le\frac{\sigma(N)}{N^{1+\epsilon}}\qquad(k\ge1).
\]
The prime factorization of the maximizer changes at the critical parameters
\begin{equation}\label{eq:critical}
 \epsilon(p,a)=\frac{\log\bigl((p^{a+1}-1)/(p^{a+1}-p)\bigr)}{\log p},
 \qquad p\text{ prime},\ a\ge1.
\end{equation}
Ordering the critical parameters in decreasing order gives the usual construction of CA numbers \cite{AE,Briggs}. At a transition involving a single prime $p$, the next maximizer is obtained by multiplication by $p$. If two distinct primes $p<q$ have the same critical parameter, the four tied maximizers are $N,pN,qN,pqN$; the middle quotient in their increasing order is $q/p$. No three distinct primes can have the same critical parameter, by the result of Siegel quoted in \cite[p.~455]{AE}. No simultaneous transition was observed in our computations. We do not assume that all critical parameters are distinct.

There is also a geometric interpretation \cite{Musin}. Place the integer $k$ at $(\log k,\log k-\log\sigma(k))$. The maximizers of $\sigma(k)/k^{1+\epsilon}$ are the points on a lower supporting line of slope $-\epsilon$. Unique maximizers are exposed vertices; tied maximizers lie on a common supporting segment.

Following Wu \cite{Wu19}, let $p(n)$ be the largest prime factor of a CA number and $q(n)$ the next prime. The three classes are
\begin{equation}\label{eq:classes}
 \begin{array}{ll}
 \CA1:&\log n\le p(n),\\
 \CA2:&p(n)<\log n<q(n),\\
 \CA3:&\log n\ge q(n).
 \end{array}
\end{equation}
\begin{theorem}[Wu \cite{Wu19}, Corollary~4]\label{th:wu}
The Riemann hypothesis is true if and only if
\[
 G(n)<\eg\qquad\text{for every }n\in\CA2\text{ with }n>5040.
\]
\end{theorem}
Wu's comparison results explain the reduction. Apart from the initial small cases, a step ending at a CA1 number decreases $G$, while a CA3 number has a smaller value of $G$ than the first subsequent CA2 number. Together with the initial range, these comparisons imply that a violation of Robin's inequality is detected on CA2. It is also enough to establish the inequality for all sufficiently large CA2 numbers: if RH is false, there are arbitrarily large CA counterexamples, and Wu's comparisons, together with the unboundedness of CA2, give arbitrarily large CA2 counterexamples.

In Wu's sample of $5\,763\,320$ CA numbers, $120\,529$ are CA1, $5\,565$ are CA2, and $5\,637\,226$ are CA3. The CA2 members therefore form less than one thousandth of this sample.

\section{Experiments and the questions they suggest}
\label{sec:experiments}
The figures compare $G$, $T$, and $R$ on CA and CA2, using $u=\log\log n$ as the horizontal coordinate. In the computations, $R$ is replaced by a finite sum over zeros; the labels $R$ and $R-T$ in the figures refer to this finite sum. The truncation error is discussed in Section~\ref{sec:computation}. Lines join the sampled values.

\subsection{Why CA2 stands out}
\begin{figure}[htbp]
\centering
\includegraphics[width=.60\linewidth]{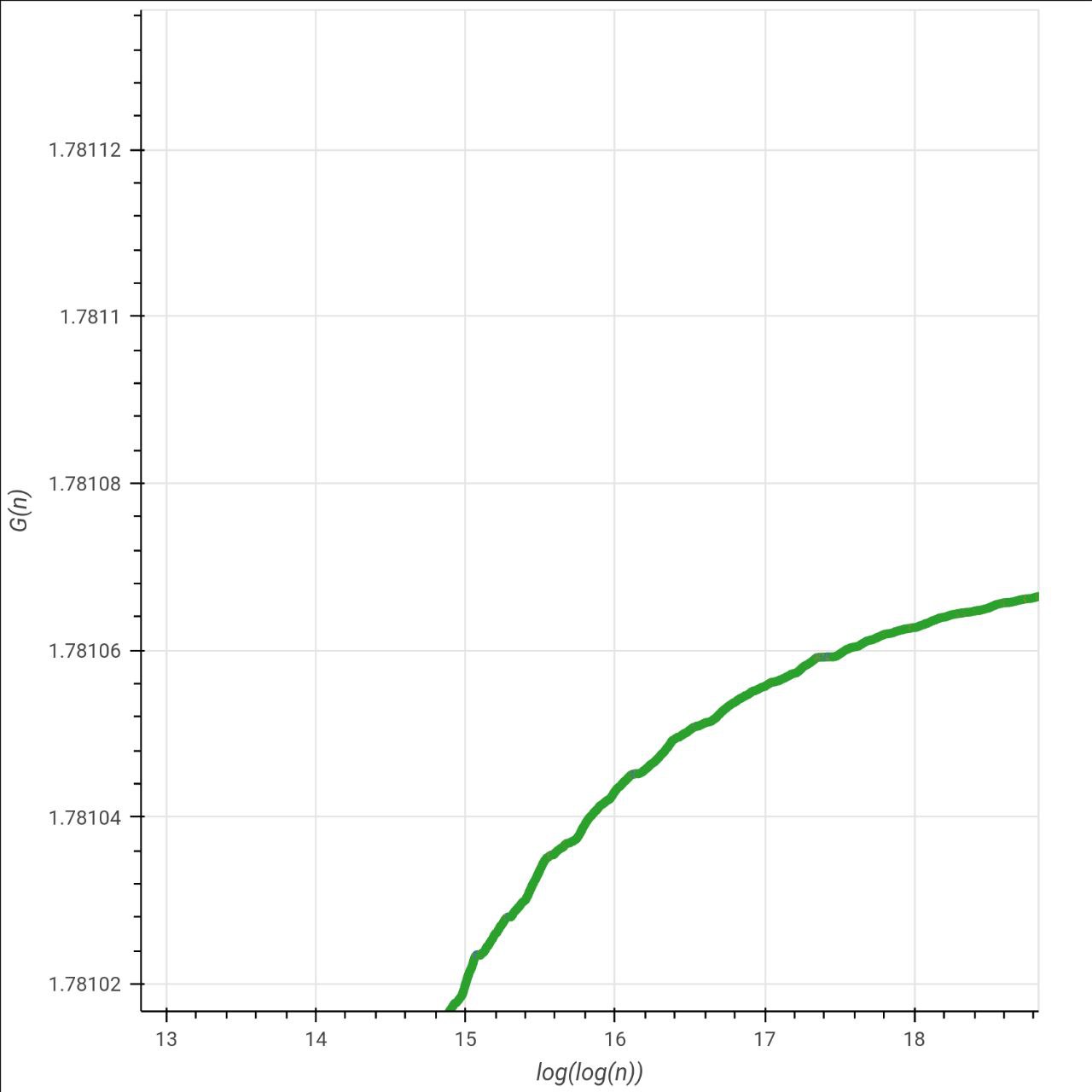}
\caption{Overview of $G(n)=\sigma(n)/(n\log\log n)$. The displayed values remain below $\eg$. On this scale the local oscillations are difficult to see and most points appear green (CA3).}
\label{fig:G-overview}
\end{figure}

\begin{figure}[htbp]
\centering
\includegraphics[width=.49\linewidth]{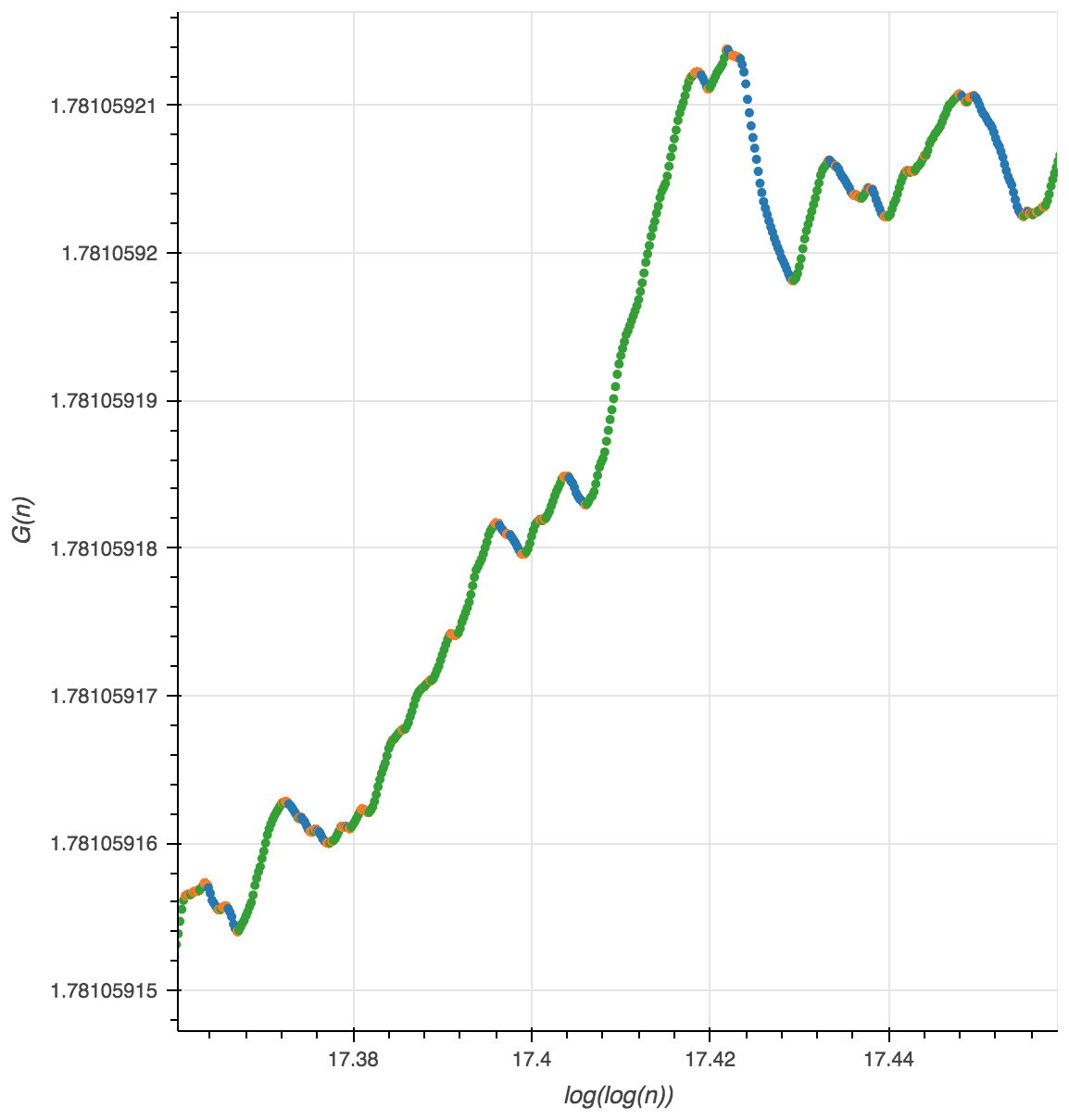}\hfill
\includegraphics[width=.49\linewidth]{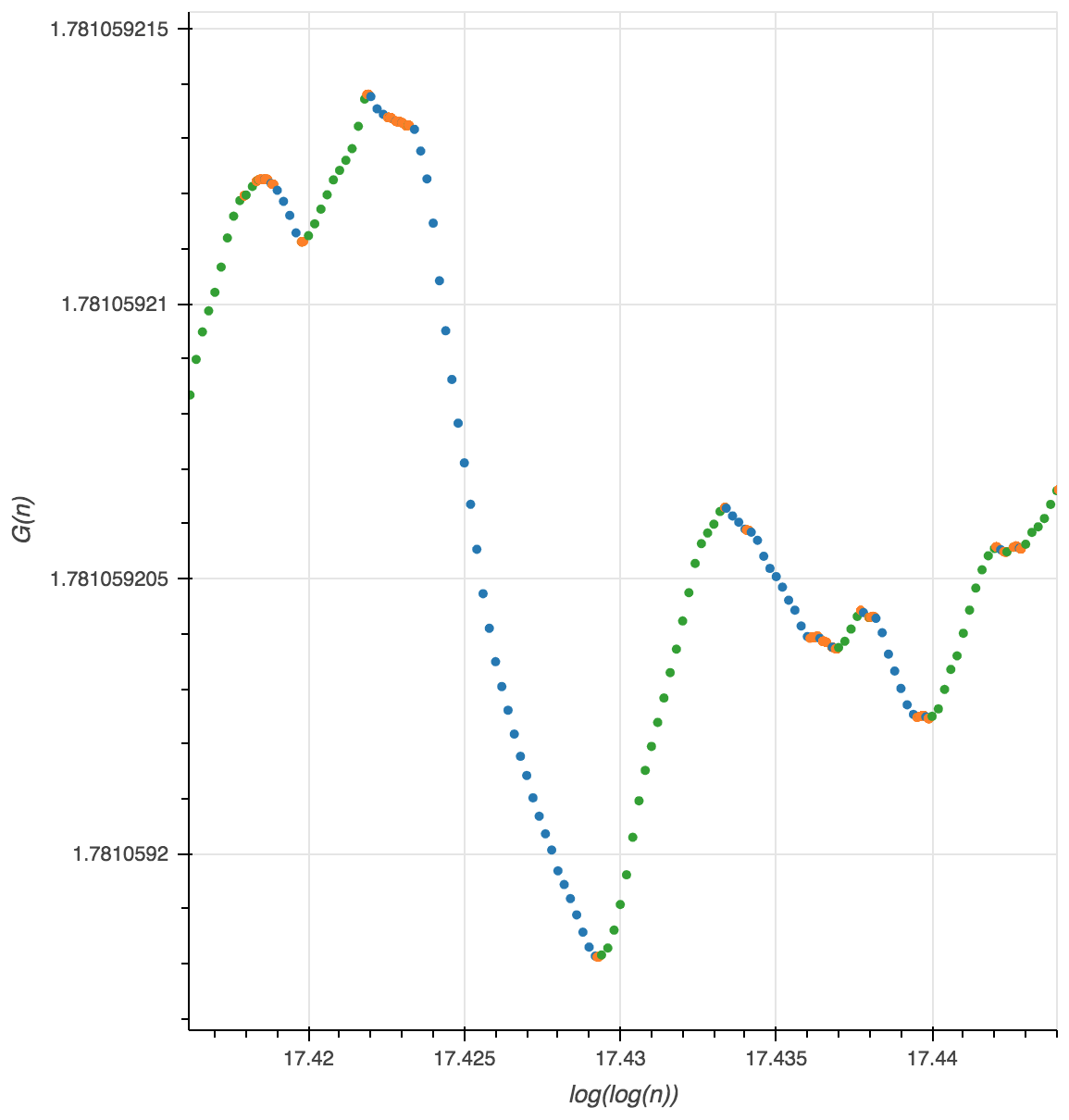}
\caption{Enlarged views of $G$. Blue CA1 points occur on descending portions, while orange/red CA2 points appear near changes of behavior. These plots resolve variations that are not visible in Figure~\ref{fig:G-overview}.}
\label{fig:G-details}
\end{figure}

The colors identify the class of each plotted point. In Figures~\ref{fig:G-overview}, \ref{fig:G-details}, and~\ref{fig:D-colored}, blue denotes CA1, orange/red denotes CA2, and green denotes CA3. The overview of $G$ shows an increasing approach towards $\eg=1.7810724\ldots$, dominated visually by CA3 points. The enlarged views reveal much smaller rises and falls within this broad trend.

These plots illustrate the comparisons behind Theorem~\ref{th:wu}: CA1 points lie on descending portions, and CA2 points occur near changes of behavior. The first subsequent CA2 number dominates a CA3 number even when the path between them is not monotone. Wu \cite[Theorem~3(1)]{Wu19} also proves that a transition by a single prime $N_{i+1}=pN_i$ starting in CA3 increases $G$ when $p\ge3$.

Wu leaves the case $p=2$ open. The plots motivate the broader question of whether every strict local maximum of $G$ belongs to CA2. The results just cited do not settle this: a CA3 number followed by multiplication by $2$, as well as simultaneous transitions by two primes, still require examination.

\begin{question}[The remaining local step]\label{q:local}
If $N_i$ is CA3 and the next CA number is $N_{i+1}=2N_i$, must $G(N_{i+1})>G(N_i)$ hold?
\end{question}

\FloatBarrier
\subsection{Can the test set be reduced further?}
\label{sec:thinning}
One may retain only the record values of $G$ on CA2 above $5040$. Write
\[
 \mathcal A=\{n\in\CA2:n>5040\},
\]
and define

\begin{equation}\label{eq:records}
 \mathcal B=\{n\in\mathcal A:G(n)>G(k)\text{ for every }k\in\mathcal A,\ k<n\}.
\end{equation}
The least member of $\mathcal A$ is included. For each $n\in\mathcal A$, the maximum of $G$ on the finite set $\mathcal A\cap[1,n]$ is attained at a record point: take its first maximizer. Hence Robin's inequality holds on $\mathcal A$ if and only if it holds on $\mathcal B$. By Theorem~\ref{th:wu}, testing the inequality on $\mathcal B$ is equivalent to RH.

This reduction is useful computationally only if the records can be located efficiently. A direct search still evaluates $G$ on every member of $\mathcal A$. Applying the same record rule to $\mathcal B$ removes no further points, so successive reductions require a different selection rule.

\begin{question}[Successive reductions]\label{q:thinning}
Can one construct explicit nested test sets
\[
 \{n\in\CA:n>5040\}\supseteq\mathcal A_1=\mathcal A
 \supseteq\mathcal A_2\supseteq\cdots
\]
such that each inclusion is strict and Robin's inequality on each fixed set remains equivalent to RH? What can be said about the intersection of these sets?
\end{question}
For a single reduction, it is sufficient that each discarded member either satisfies Robin's inequality or has $G$ no larger than that of a retained member.

The higher-order CA sets $C^{(k)}(\mathcal T)$ in \cite{Musin26} provide related constructions on the full CA sequence. They are defined by supporting lines after changes of coordinates. Under the stated conditions, the sets are nested and Robin's inequality on any fixed level is equivalent to RH. Some of these families have infinitely many members at every level and empty intersection, unconditionally. Their construction does not assert that they lie in Wu's CA2 class. Thus the relation between these families and successive reductions within CA2 remains to be studied. An empty intersection alone would not prove RH: counterexamples could persist at every level while each fixed integer is eventually removed.

\subsection{Oscillation and a slowly varying difference}
\begin{figure}[htbp]
\centering
\includegraphics[width=.80\linewidth]{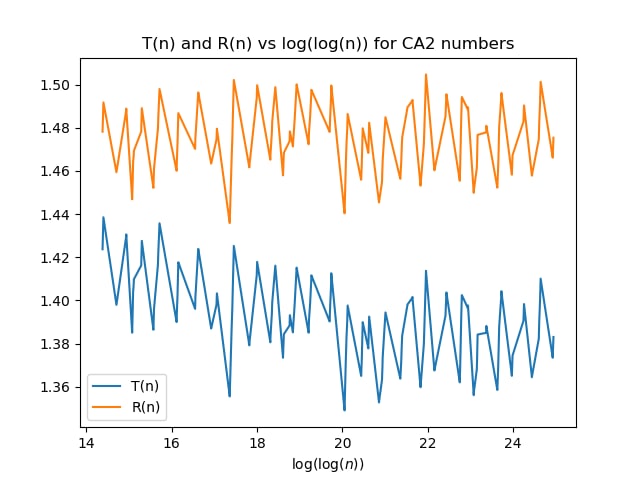}
\caption{The curves labeled $T(n)$ and $R(n)$ on the CA2 sample, with horizontal coordinate $u=\log\log n$. Their visible oscillations are similar, but the vertical separation changes across the range. The figure includes every point of our CA2 sample.}
\label{fig:TR}
\end{figure}

Figure~\ref{fig:TR} compares $T$ with the sum over zeros on CA2. Both curves display a repeated oscillatory pattern across the sample, while the upper curve is separated from the lower one by a positive gap. This separation is easier to study directly than through the two curves themselves.

\begin{figure}[htbp]
\centering
\includegraphics[width=.49\linewidth]{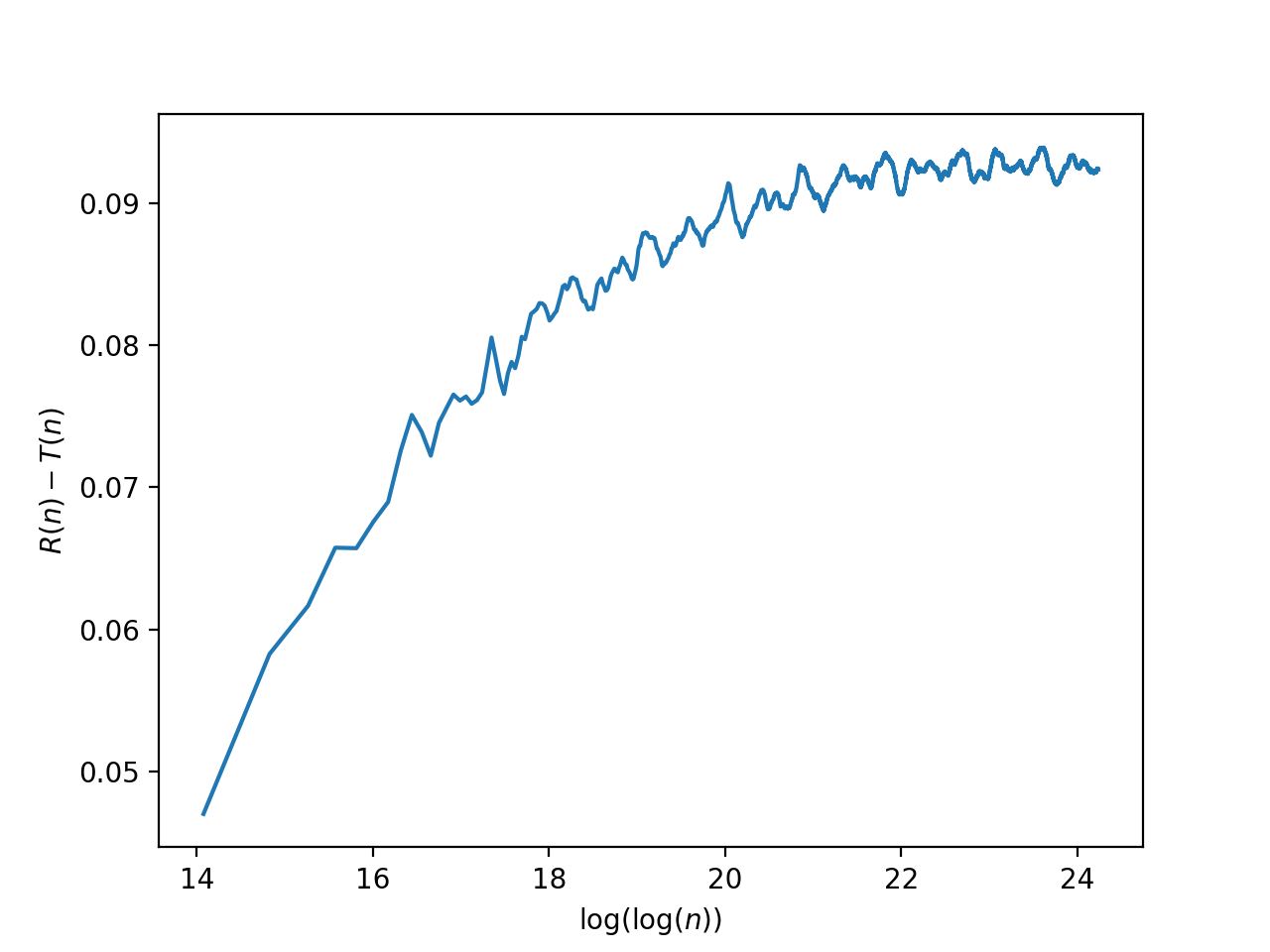}\hfill
\includegraphics[width=.49\linewidth]{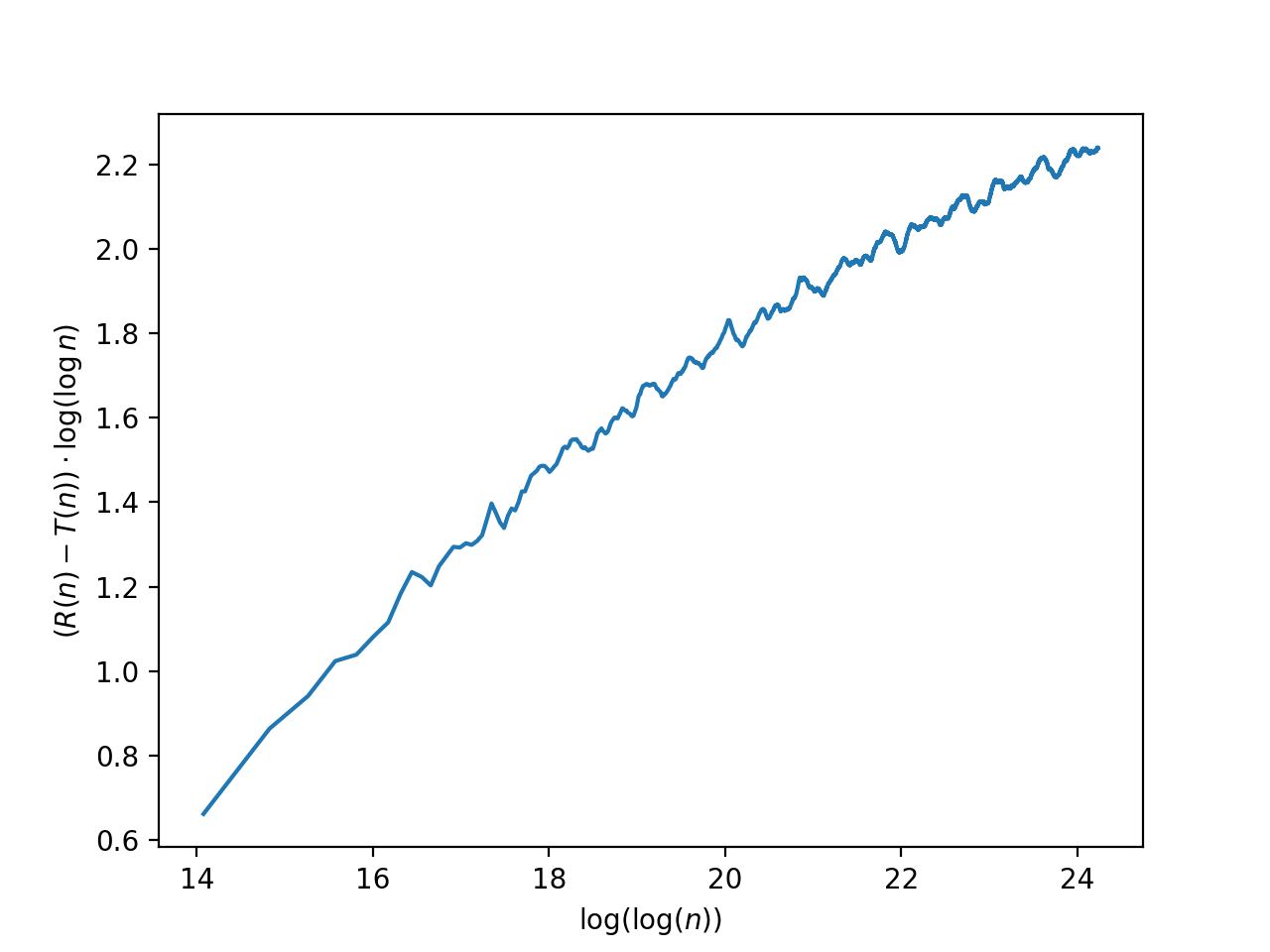}
\caption{Residual scales on sampled CA numbers. Left: $D=R-T$, sampled every $10\,000$ CA numbers. Right: $uD=\eg Q$, sampled every $1000$ CA numbers. The vertical scales differ by the factor $u$ and the two panels use different sampling rates.}
\label{fig:residual-scale}
\end{figure}

For the sampled full CA sequence, the left panel of Figure~\ref{fig:residual-scale} shows the residual increasing from a little below $0.05$ to a little above $0.09$ over $14\lesssim u\lesssim25$, with smaller local oscillations. The increase becomes less pronounced near the right end of the plot. The right panel shows the same type of experiment on the scale $uD(n)$: its values rise from roughly $0.65$ to $2.2$. These values are read approximately from the displayed curves.

After division by $\eg$, the right panel corresponds approximately to $0.36\lesssim Q\lesssim1.25$ on its displayed sample. The main experimental feature is the drift across this interval. A nearly constant value of $D$ over a short interval does not imply that $Q$ is stabilizing; multiplication by $u$ can still reveal substantial variation. Accordingly, the finite-range increase does not rule out eventual boundedness or oscillation.

\begin{question}[The residual coefficient]\label{q:limit}
Assume RH. Does $Q(n)$ converge as $n\to\infty$ through CA numbers? If it does not, what are its lower and upper limits? Are the answers the same on CA2?
\end{question}
If the limit is $C$, its meaning is the expansion
\[
 T(n)=R(n)-\frac{\eg C+o(1)}{\log\log n}.
\]
The displayed values do not give a stable estimate of $C$. Comparisons on successive $u$-intervals may help distinguish convergence from persistent oscillation.

One can ask for an upper bound without assuming convergence.
\begin{question}[A proposed upper bound]\label{q:bound}
Under RH, is it true that $Q(n)\le3/2$ for all sufficiently large CA numbers? Equivalently, does
\begin{equation}\label{eq:proposed-bound}
 T(n)\ge R(n)-\frac{3\eg}{2\log\log n}
\end{equation}
eventually hold along the CA sequence?
\end{question}
The plotted sample is consistent with this bound after the normalization $Q=e^{-\gamma}uD$. It does not determine an eventual threshold or the best coefficient.

\subsection{What is retained on CA2?}
\begin{figure}[htbp]
\centering
\includegraphics[width=.49\linewidth]{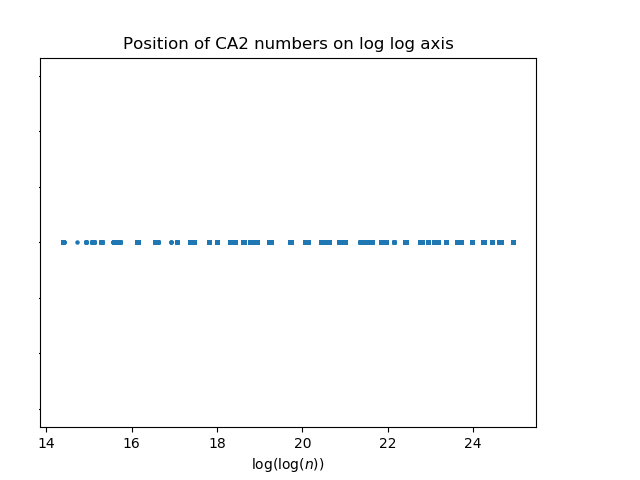}\hfill
\includegraphics[width=.49\linewidth]{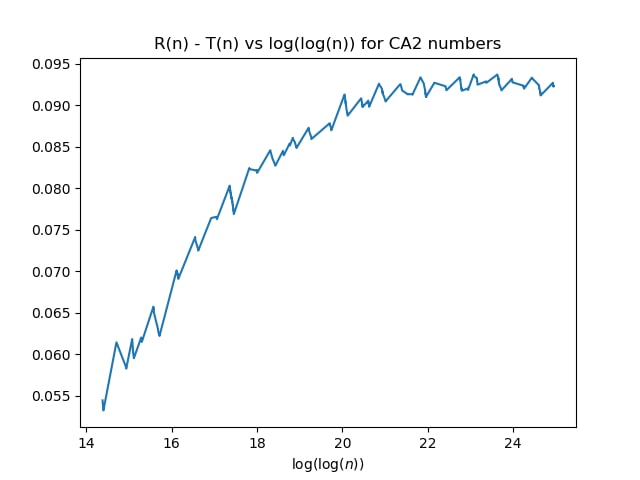}
\caption{The CA2 sample. Left: positions on the $u$ axis, with the meaningless vertical scale omitted. Right: the residual labeled $R-T$ on this sample. Its broad trend resembles the full CA residual, while the irregular spacing affects the appearance of the local oscillations.}
\label{fig:CA2}
\end{figure}

Figure~\ref{fig:CA2} shows the positions of the CA2 sample and its residual. The points are sparse and unevenly spaced on the $u$ axis. Nevertheless, the broad increase and the small oscillations of the residual resemble those of the entire CA sequence sampled in Figure~\ref{fig:residual-scale}.

The comparison suggests using CA2 to study some features of the residual at a much lower sampling cost. However, there is a distinction between the known reduction for $G$ and the proposed use of CA2 for $Q$. Robin's inequality depends on whether $G$ crosses the constant level $\eg$; the residual also subtracts the oscillatory function $R$ and changes the scale. Wu's reduction therefore does not place the extrema of $Q$ in CA2.

\begin{question}[CA2 and the residual envelope]\label{q:CA2}
Under RH, do the lower and upper limits of $Q$ on CA2 agree with the corresponding limits on all CA numbers?
\end{question}
A direct numerical test would compare both samples in identical $u$-windows and record their minimum, maximum, and number of points. The existing curves suggest the comparison, but do not supply those window statistics.
\FloatBarrier

\subsection{Local behavior of the residual}
\begin{figure}[htbp]
\centering
\includegraphics[width=.49\linewidth]{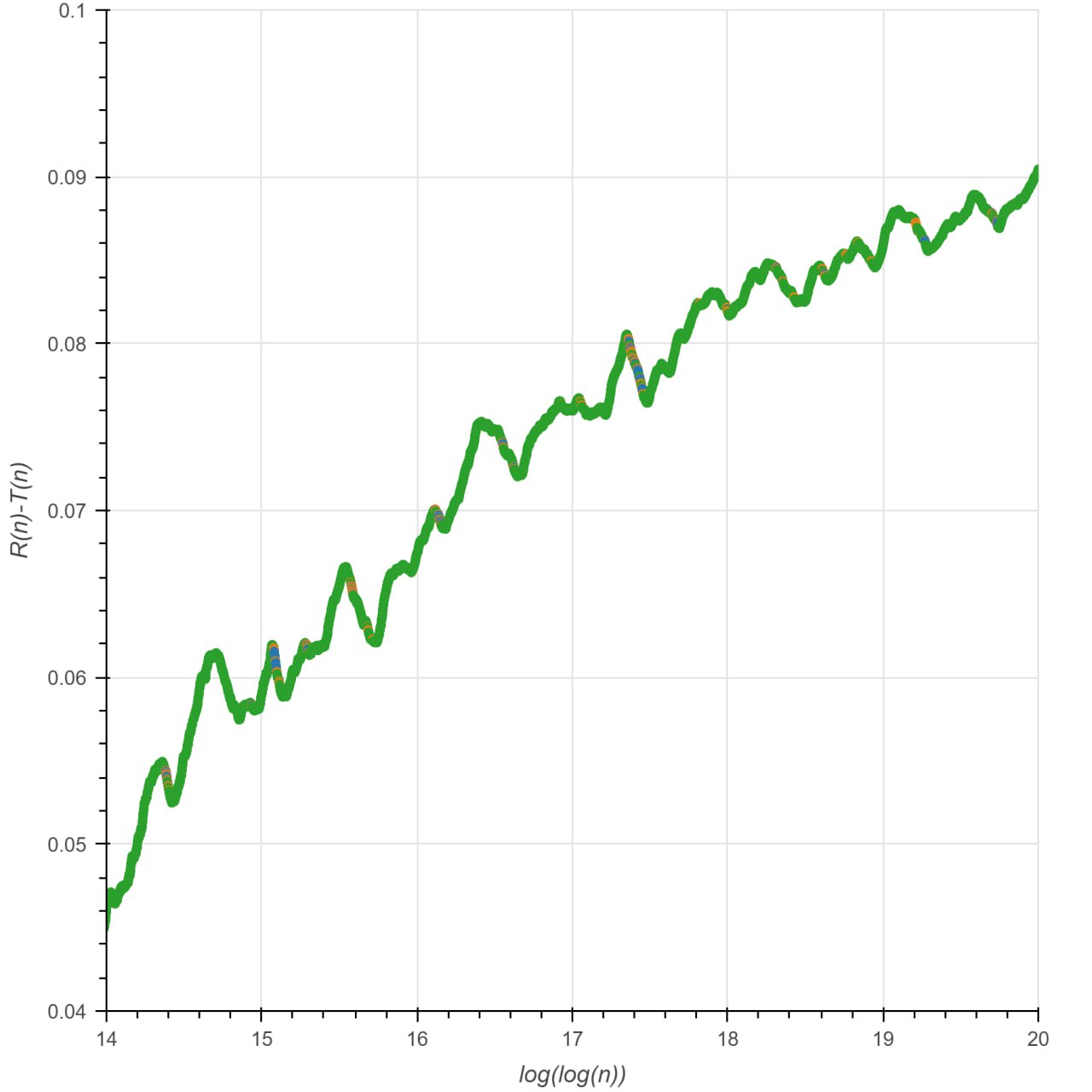}\hfill
\includegraphics[width=.49\linewidth]{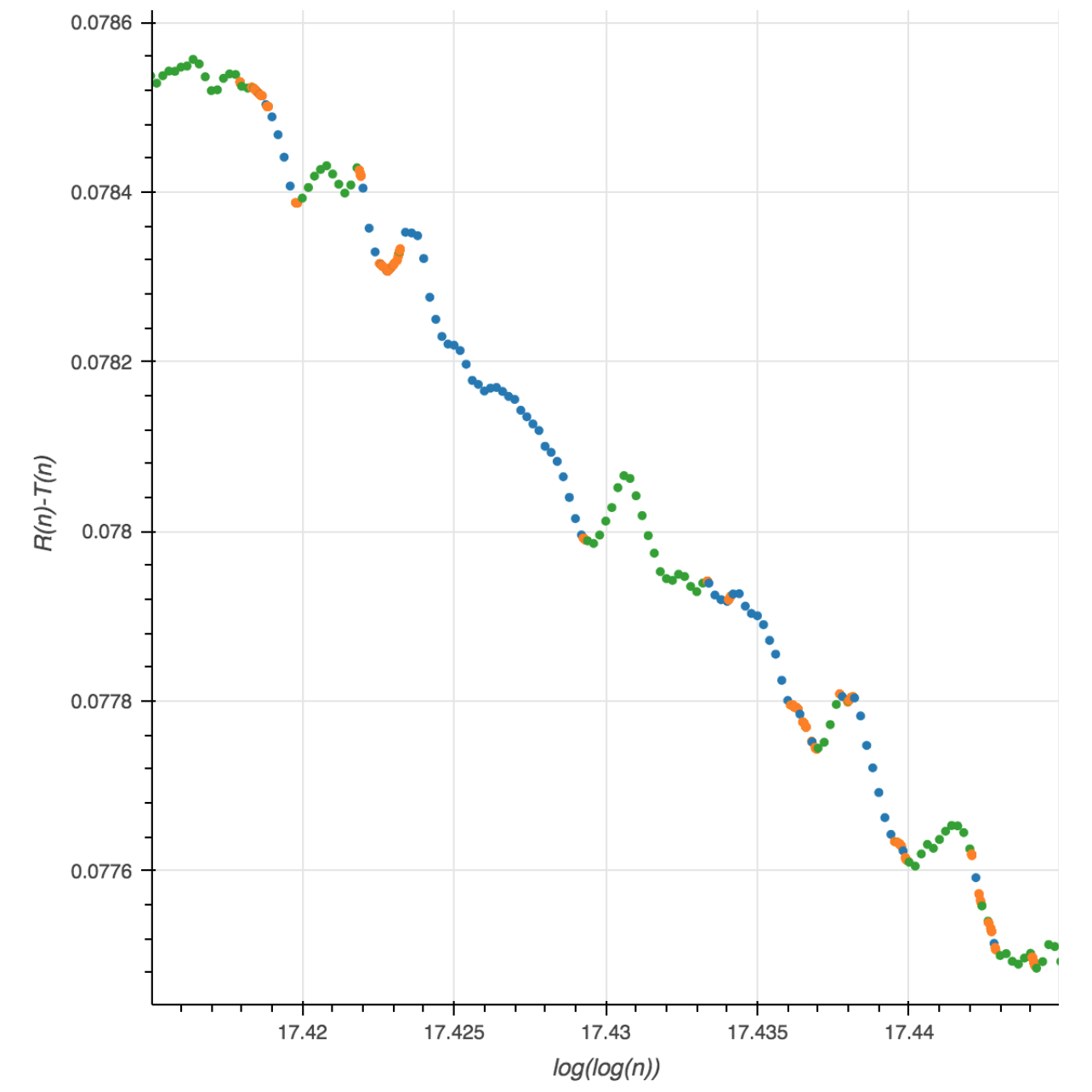}
\caption{Class-colored residual: overview and detail. The vertical axes show $D=R-T$. The normalization $P=e^{-\gamma}D$ preserves the signs of consecutive differences.}
\label{fig:D-colored}
\end{figure}

For the residual, the class pattern is more complicated. Figure~\ref{fig:D-colored} juxtaposes an increasing trend with a detail in which blue and green points lie on descending portions. Thus, a class rule for increases of $G$ cannot simply be transferred to $D$. Multiplication by the positive constant $e^{-\gamma}$ preserves this local pattern for $P$; the additional factor $u$ in $Q$ may change it.

\FloatBarrier
These observations suggest examining consecutive differences conditional on both the CA class and the transition prime. Such a comparison would separate a systematic class effect from the position within an oscillation of the sum over zeros. The figures show examples of this local behavior; a systematic comparison requires the consecutive numerical values.
\FloatBarrier

\section{Computational outline and accuracy}
\label{sec:computation}
\subsection{Generation and sampling}
The computational method follows the critical-parameter description of CA numbers. If a prime currently has exponent $a$, its next event is $\epsilon(p,a+1)$ from \eqref{eq:critical}; the largest remaining critical parameter is processed next. Storing the prime exponents, $L=\log n$, and $A=\sigma(n)/n$ avoids storing the extremely large integer $n$. At a transition $n\mapsto pn$, with current exponent $a\ge0$, the updates are
\[
 L\longleftarrow L+\log p,\qquad
 A\longleftarrow A\frac{1-p^{-(a+2)}}{1-p^{-(a+1)}}.
\]
Then $T=(\eg\log L-A)\sqrt L$. Comparisons of nearly equal event parameters and the cancellation in the evaluation of $T$ were the two main precision issues.

The computations used C++ and Python programs, with 100-bit Arb arithmetic in the C++ generator; its first million transition factors were checked against an independent generator. Transition factors were generated through CA index $10^{10}$, and $T$ was evaluated through index $24\cdot10^6$. Generation of the transition factors and evaluation of $T$ are separate computations with different ranges. The sampling rates for the full CA sequence are given in the figure captions; the CA2 figures use all points in our CA2 sample.

\subsection{Truncation of the sum over zeros}
Let $R_K$ denote the sum in \eqref{eq:R} truncated after the first $K$ positive ordinates. Under RH, the elementary tail estimate is
\begin{equation}\label{eq:tail}
 \begin{split}
 E_K&=\tau-2\sum_{j=1}^K\frac1{\tfrac14+t_j^2},\\
 |R(n)-R_K(n)|&\le\eg E_K,\qquad
 |Q(n)-Q_K(n)|\le uE_K,
 \end{split}
\end{equation}
where $Q_K=e^{-\gamma}(R_K-T)u$. This follows by bounding each omitted cosine by one. The bound for $R$ is uniform in $n$, whereas the bound for $Q$ grows with $u$.

For example, evaluation using Odlyzko's ordinate tables \cite{Zeros} gives the rounded scales in Table~\ref{tab:tail}. The last column applies to $u\le25$, which includes the total interval displayed.
\begin{table}[htbp]
\centering
\begin{tabular}{rrrr}
\toprule
$K$ & $E_K$ & $\eg E_K$ & $25E_K$\\
\midrule
$1\,000$   & $1.43916\cdot10^{-3}$ & $2.56324\cdot10^{-3}$ & $3.59789\cdot10^{-2}$\\
$10\,000$  & $2.69398\cdot10^{-4}$ & $4.79818\cdot10^{-4}$ & $6.73496\cdot10^{-3}$\\
$100\,000$ & $4.41272\cdot10^{-5}$ & $7.85938\cdot10^{-5}$ & $1.10318\cdot10^{-3}$\\
\bottomrule
\end{tabular}
\caption{Approximate zero-tail scales under RH. The third column concerns $R$; the fourth concerns $Q$ when $u\le25$. These rounded evaluations do not include ordinate rounding, arithmetic errors, or an error in $T$.}
\label{tab:tail}
\end{table}

We chose $K=1000$ when performing our computations of the sum over zeros. At that value, the tail scale alone for $Q$ near $u=25$ is about $0.036$. Increasing to $K=10\,000$ or $100\,000$ reduces this scale to approximately $0.00673$ or $0.00110$, respectively. These are the scales for $Q$ in the last column, distinct from the smaller errors for $R$. Resolving stabilization, or the higher-order coefficients, requires a tail scale below the changes being measured, together with control of arithmetic errors. In particular, a fixed $K$ does not give a uniform error bound for $Q$ as $u\to\infty$.

The Riemann--von Mangoldt formula gives $t_j\sim2\pi j/\log j$. Consequently,
\[
 E_K=2\sum_{j>K}\frac1{\tfrac14+t_j^2}
 \sim\frac{(\log K)^2}{2\pi^2K}.
\]
Thus $K$ of order $u_H(\log u_H)^2$ keeps the bound $u_HE_K$ of constant order as the upper endpoint $u_H$ grows. To make this bound tend to zero, one may choose $K$ so that $u_H(\log K)^2/K\to0$.

\subsection{Numerical materials}
The accompanying files contain the figures used in this paper and four additional plots from earlier computations. The underlying point tables and programs are not included. The numerical ranges quoted from the figures are approximate. The additional plots are not used here because their axis conventions or normalizations differ from those used in the main text.

\section{Successive corrections and asymptotic questions}
\label{sec:asymptotic}
\subsection{From Nicolas's refinement to the next correction}
\label{sec:nicolas}
Nicolas \cite{Nic21} gives an effective form of Ramanujan's conditional upper estimate for $\sigma(n)/n$. In the notation of \eqref{eq:residuals}, its asymptotic form \eqref{eq:one-sided} implies that there are constants $C>0$ and $n_0$ such that, under RH,
\begin{equation}\label{eq:N-Q}
 P(n)\le\frac{C}{u},\qquad Q(n)\le C\qquad(n\ge n_0).
\end{equation}
This motivates the scale $1/u$ for the first correction. A constant in a one-sided inequality need not be the coefficient of an asymptotic expansion. Question~\ref{q:bound} asks whether $C=3/2$ is an eventual upper bound on CA. A separate question is whether the next terms have constant coefficients or retain an oscillatory dependence on $u$.

\begin{question}[Higher-order terms]\label{q:higher}
Assume RH. Do there exist constants $A_1,A_2,\ldots$ such that, for each fixed $r\ge1$,
\begin{equation}\label{eq:expansion}
 P(n)=\frac{A_1}{u}+\frac{A_2}{u^2}+\cdots+\frac{A_r}{u^r}
       +o(u^{-r}),\qquad u=\log\log n\to\infty,
\end{equation}
along CA numbers, or along CA2? If constant coefficients are inadequate, what functions describe the next terms?
\end{question}
The associated experiment is sequential. For a candidate $A_1$, inspect $u(Q-A_1)$; if this stabilizes near $A_2$, inspect $u^2(Q-A_1-A_2/u)$. In general the quantity to examine is
\begin{equation}\label{eq:successive-residual}
 Q_r(n)=u^r\left(P(n)-\sum_{j=1}^{r-1}\frac{A_j}{u^j}\right).
\end{equation}
Thus $Q_1=Q$, and \eqref{eq:expansion} would give $Q_r\to A_r$. The present values of $Q_1$ do not give a stable first coefficient, so the higher-order differences cannot yet be interpreted numerically.

If such an expansion exists, it refines the divisor-sum formula to
\begin{equation}\label{eq:sigma-expansion}
 \frac{\sigma(n)}n=\eg\log\log n-\frac{R(n)}{\sqrt{\log n}}
 +\frac{\eg}{\sqrt{\log n}}
  \left(\sum_{j=1}^r\frac{A_j}{u^j}+o(u^{-r})\right).
\end{equation}
One-sided estimates for the remainder would give corresponding inequalities for $\sigma(n)/n$. An unconditional proof that $T(n)$ is eventually bounded below by a positive constant on CA would imply RH. For example, the estimate $T(n)\ge c_1-O(1/u)$, if proved without assuming RH, would suffice. In contrast, the definition of $R$ in \eqref{eq:R} and the comparisons above are made under RH.

The accuracy requirement also increases with each stage. By \eqref{eq:tail}, the zero-tail contribution to $Q_r$ is at most $u^rE_K$. An error in a fitted coefficient $A_j$ contributes that error multiplied by $u^{r-j}$. These amplifications should be kept below the variation under investigation before a new apparent term is interpreted mathematically.

\subsection{Sharpness of Ramanujan's constants}
The residual and the oscillatory envelope lead to different asymptotic problems. Questions~\ref{q:limit}--\ref{q:bound} concern the scale on which $T$ approaches $R$. Even if this remainder were completely understood, the extremes in behavior of the sum over zeros along CA numbers would still have to be examined to determine the exact constants in \eqref{eq:ramanujan}.

\begin{conjecture}[Sharpness of Ramanujan's constants]\label{conj:sharpness}
Under RH, the limits along CA numbers satisfy
\[
 R_L=c_1,\qquad R_U=c_2.
\]
\end{conjecture}
The cosine bound giving $c_1$ and $c_2$ does not show that these endpoints are approached along the CA sequence. A possible explanation would combine $T(n)-R(n)\to0$ along CA with simultaneous approximation of the phases $t_j\log\log n$. Linear independence over $\mathbb Q$ of the positive ordinates would allow any finite collection of phases to approach prescribed values as the real parameter $u$ varies. To apply this to CA numbers, one must also pass from real values of $u$ to the values $u=\log\log n$ with $n\in\CA$. This describes a possible route to the conjecture; these additional assertions are not consequences of the numerical plots.

For the residual problem, a useful next experiment is to divide an extended range into comparable $u$-windows, keeping the zero-tail and arithmetic errors below a prescribed tolerance for $Q$. Recording the minimum and maximum in each window would directly test the proposed upper bound; recording the center and spread would help distinguish convergence from persistent oscillation or continued drift. Using the same windows for CA and CA2 would address Question~\ref{q:CA2} without relying on the visual similarity of differently spaced samples.

For the local problem, a systematic computation of CA transitions should test the multiplier-2 steps of Question~\ref{q:local} and tabulate consecutive changes in $D$ and $Q$ by CA class and transition prime.

\bigskip
\noindent\begin{minipage}{\linewidth}
\small
Dmitry R. Leshchiner\\
\textit{email}: \href{mailto:leshch6@gmail.com}{leshch6@gmail.com}

\medskip
\noindent Oleg R. Musin\\
University of Texas Rio Grande Valley,\\
Brownsville, TX 78520, USA.\\
\textit{email}: \href{mailto:oleg.musin@utrgv.edu}{oleg.musin@utrgv.edu}

\medskip
\noindent James Thompson\\
\textit{email}: \href{mailto:james.thompson.5548@gmail.com}{james.thompson.5548@gmail.com}
\end{minipage}
\end{document}